\input amstex
\documentstyle{amsppt}
\magnification=\magstep1                        
\hsize6.5truein\vsize8.9truein                  
\NoRunningHeads
\loadeusm

\magnification=\magstep1                        
\hsize6.5truein\vsize8.9truein                  
\NoRunningHeads
\loadeusm

\document
\topmatter

\title
The sequence of partial sums of a unimodular power series is not ultraflat 
\endtitle

\rightheadtext{ultraflat unimodular polynomials}

\author Tam\'as Erd\'elyi
\endauthor

\address Department of Mathematics, Texas A\&M University,
College Station, Texas 77843, College Station, Texas 77843 \endaddress

\endthanks


\date April 23, 2025
\enddate

\email terdelyi\@tamu.edu
\endemail

\dedicatory \enddedicatory

\abstract
We show that if $(a_j)_{j=0}^\infty$ is a sequence of numbers $a_j \in {\Bbb C}$ with $|a_j|=1$, and
$$P_n(z) = \sum_{j=0}^n{a_jz^j}\,, \qquad n=0,1,2,\ldots\,,$$
then $(P_n)$ is NOT an ultraflat sequence of unimodular polynomials.
This answers a question raised by Zachary Chase.
\endabstract

\endtopmatter

\head 1. Introduction \endhead

Let
$${\Cal K}_n := \bigg\{Q_n: Q_n(z) = \sum_{k=0}^n{a_k z^k}, \enskip  a_k \in {\Bbb C}\,,\enskip |a_k| = 1 \bigg\}\,.$$
The class ${\Cal K}_n$ is called the collection of all complex unimodular polynomials of degree $n$.
Elements of ${\Cal K}_n$ may be called Kahane polynomials of degree $n$.  
Let
$${\Cal L}_n := \bigg\{Q_n: Q_n(z) = \sum_{k=0}^n{a_k z^k}, \enskip  a_k \in \{-1,1\} \bigg\}\,.$$
The class ${\Cal L}_n$ is called the collection of all real unimodular polynomials of degree $n$.
Elements of ${\Cal L}_n$ are called Littlewood polynomials of degree $n$.
By Parseval's formula,
$$\int_{0}^{2\pi}{\left| P_n(e^{it}) \right|^2 \, dt} = 2\pi(n+1)$$
for all $P_n \in {\Cal K}_n$. Therefore
$$\min_{t \in [0,2\pi]}{\left|P_n(e^{it})\right|} \leq \sqrt{n+1} \leq \max_{t \in [0,2\pi]}{\left|P_n(e^{it})\right|}\,.$$
An old problem (or rather an old theme) is the following.

\proclaim{Problem 1.1 (Littlewood's Flatness Problem)}
How close can a polynomial $P_n \in {\Cal K}_n$ or $P_n \in {\Cal L}_n$ come to satisfying
$$\left|P_n(e^{it})\right| = \sqrt{n+1}\,, \qquad t \in {\Bbb R}? \tag 1.1$$
\endproclaim

Obviously (1.1) is impossible if $n \geq 1$. So one must look for less than (1.1), but then there are 
various ways of seeking such an ``approximate situation". One way is the following.
In his paper [Li1] Littlewood had suggested that, conceivably, there might exist a sequence
$(P_n)$ of polynomials $P_n \in {\Cal K}_n$ (possibly even $P_n \in {\Cal L}_n$) such
that $(n+1)^{-1/2}|P_n(e^{it})|$ converge to $1$ uniformly in $t \in {\Bbb R}$.
We shall call such sequences of unimodular polynomials ``ultraflat". More precisely, we give the
following definition.

\proclaim{Definition 1.2} Given a positive number $\varepsilon$, we say that a polynomial
$P_n \in {\Cal K}_n$ is $\varepsilon$-flat if
$$(1-\varepsilon)\sqrt{n+1} \leq \left|P_n(e^{it})\right| \leq (1+\varepsilon)\sqrt{n+1}\,, \qquad t \in {\Bbb R}\,.$$
\endproclaim

\proclaim{Definition 1.3} Given a sequence $(\varepsilon_{n})$ of positive numbers tending to $0$, we say that a sequence
$(P_n)$ of polynomials $P_n \in {\Cal K}_n$ is $(\varepsilon_n)$-ultraflat if each $P_n$ is $(\varepsilon_n)$-flat.
We simply say that a sequence $(P_n)$ of polynomials $P_n\in {\Cal K}_n$ is
ultraflat if it is $(\varepsilon_n)$-ultraflat with a  sequence $(\varepsilon_n)$ of positive numbers converging to $0$.
\endproclaim

The existence of an ultraflat sequence of unimodular polynomials seemed very unlikely, in view of a
1957 conjecture of P. Erd\H os (Problem 22 in [Er]) asserting that, for all
$P_n \in {\Cal K}_n$ with $n \geq 1$,
$$\max_{t \in {\Bbb R}}{\left|P_n(e^{it})\right|} \geq (1 + \varepsilon) \sqrt{n+1}\,, \tag 1.2$$
where $\varepsilon > 0$ is an absolute constant (independent of $n$).
Yet, refining a method of K\"orner [K\"o], Kahane [Ka] proved that there exists
a sequence $(P_n)$ with $P_n \in {\Cal K}_n$ which is $(\varepsilon_n)$-ultraflat, where
$\varepsilon_n = O\left(n^{-1/17} \sqrt{\log n} \right)\,.$
(Kahane's paper contained though a slight error which was corrected in [QS2].)
Thus the Erd\H os conjecture (1.2) was disproved for the classes ${\Cal K}_n$.
For the more restricted class ${\Cal L}_n$ the analogous Erd\H os conjecture
is unsettled to this date. It is a common belief that the analogous Erd\H os conjecture
for ${\Cal L}_n$ is true, and consequently there is no ultraflat sequence of polynomials 
$P_n \in {\Cal L}_n$.
For an account of some of the work done till the mid 1960's, see Littlewood's book [Li2]
and [QS2]. Properties of the Rudin-Shapiro polynomials have played a central role in [BBM] 
as well as in [Er8] to prove a longstanding conjecture of Littlewood on the existence of flat
Littlewood polynomials $S_n \in {\Cal L}_n$ satisfying the inequalities
$$c_1 \sqrt{n+1} \leq \left|S_n(e^{it})\right| \leq c_2 \sqrt{n+1}\,, \qquad t \in {\Bbb R}\,, \quad n=0,1,2,\ldots \,,$$
with absolute constants $c_1 > 0$ and $c_2 > 0$. 
The papers [BB], [Er1]--[Er8], [EN], [Ka], [Od], [QS1], [QS2], and [Sa] deal with topics closely related 
to ultraflat sequences of unimodular polynomials.

\head 2. New Result \endhead

\proclaim{Theorem 2.1}
If $(a_j)_{j=0}^\infty$ is a sequence of numbers $a_j \in {\Bbb C}$ with $|a_j|=1$, and
$$P_n(z) = \sum_{j=0}^n{a_jz^j}\,, \qquad n=0,1,2,\ldots\,,$$
then $(P_n)$ is NOT an ultraflat sequence of unimodular polynomials.
\endproclaim

In other words, the coefficients of ultraflat unimodular polynomials must vary with the degree.
This answers a question raised by Zachary Chase.

\head Proof of Theorem 2.1 \endhead
To prove Theorem 2.1 we need the second equality of the following lemma proved in [Er1]. 

\proclaim{Lemma 2.2 (The Bernstein Factors)}
Let $q$ be an arbitrary positive real number. Let $(P_n)$ be a fixed ultraflat sequence of polynomials $P_n \in {\Cal K}_n$.
We have
$$\frac{\int_0^{2\pi}{\left| P_n^\prime(e^{it}) \right|^q \, dt}}
{\int_0^{2\pi}{ \left| P_n(e^{it}) \right|^q \, dt}}  = \frac{n^{q}}{q+1} + o_{n,q}n^{q}\,,$$
with suitable constants $o_{n,q}$ converging to $0$ as $n$ tends to $\infty$ for every fixed $q > 0$,
and 
$$\frac{\max_{0 \leq t \leq 2\pi}{\left| P_n^\prime(e^{it}) \right|}}{\max_{0 \leq t \leq 2\pi}{\left|P_n(e^{it})\right|}} = n + o_n n$$
with suitable constants $o_n$ converging to $0$ as $n$ tends to $\infty$.
\endproclaim 

\demo{Proof of Theorem 2.1}
Suppose to the contrary that $(a_j)_{j=0}^\infty$ is a sequence of numbers $a_j \in {\Bbb C}$ with $|a_j|=1$,  
$$P_n(z) = \sum_{j=0}^n{a_jz^j}\,, \qquad n=0,1,2,\ldots\,,$$
and $(P_n)$ is an $(\varepsilon_n)$-ultraflat sequence of unimodular polynomials $P_n \in {\Cal K}_n$.
Associated with $P_n \in {\Cal K}_n$ we define $Q_0(z) \equiv 1$ and 
$$Q_{n+1}(z) := 1 + z^{n+1}P_n(1/z) = 1 + \sum_{j=1}^{n+1}{a_{n+1-j}z^j}\,, \qquad n=0,1,2,\ldots\,.$$
Observe that $(Q_n)$ is an $(\varepsilon_n^*)$-ultraflat sequence of unimodular polynomials $Q_n \in {\Cal K}_n$ with 
positive real numbers $\varepsilon_n^*$ converging to $0$ and  
$$z^nQ_{n+1}^\prime(1/z) = z^n\sum_{j=1}^{n+1}{a_{n+1-j}j(1/z)^{j-1}} = \sum_{j=1}^{n+1}{a_{n+1-j}jz^{n+1-j}} = 
\sum_{k=0}^n{a_k(n+1-k)z^k}\,,$$
hence
$$\sum_{k=0}^n{P_k(z)} = \sum_{k=0}^n{(n+1-k)a_kz^k} = z^nQ_{n+1}^{\prime}(1/z)\,, \qquad 0 \neq z \in {\Bbb C}\,. \tag 2.1$$
Using Definition 1.3 we have 
$$\split \bigg| \sum_{k=0}^n{P_k(e^{it})} \bigg| & \leq \sum_{k=0}^n{\left|P_k(e^{it})\right|} \leq \, \sum_{k=0}^n{(1+\varepsilon_k)\sqrt{k+1}} \cr 
& \leq  \int_{1}^{n+2}{x^{1/2} \, dx} + \sum_{k=0}^n{\varepsilon_k}(n+1)^{1/2} \cr 
& \leq \frac 23 \, (1+\varepsilon_n^{**})(n+2)^{3/2}\,, \qquad t \in {\Bbb R}\,, \qquad n=0,1,2,\ldots\,,\cr \endsplit \tag 2.2$$ 
where 
$$\varepsilon_n^{**} := \frac{3}{2(n+2)} \sum_{k=0}^n{\varepsilon_k}\,, \qquad n=0,1,2,\ldots\,,$$
and hence the sequence $(\varepsilon_n^{**})$ of positive real numbers converges to $0$. 
On the other hand, combining (2.1) and the second equality of Lemma 2.2 applied to the ultraflat sequence 
$(Q_n)$ of unimodular polynomials $Q_n \in {\Cal K}_n$, we obtain that there are $t_n \in [0,2\pi)$ such that  
$$\bigg| \sum_{k=0}^n{P_k(e^{it_n})} \bigg| = \left| Q_{n+1}^{\prime}(e^{-it_n}) \right| \geq \frac 34 \, (n+2)^{3/2} \tag 2.3$$
for every sufficiently large $n$. Observe that (2.3) contradicts (2.2).
\qed \enddemo

\enddocument